\documentclass[fleqn,10pt]{article}
\usepackage{amsmath}
\usepackage{latexsym}
\usepackage{rotate,epsfig}
\usepackage{lscape}
\usepackage{bm}
\usepackage{graphicx}
\usepackage{color}
\newtheorem{theorem}{\bf Theorem}[section]

\newtheorem{remark}[theorem]{\bf Remark}

\renewcommand{\vartheta}{\tau}

\begin{document}
\title{Weibull analysis with sequential order statistics under a power trend model for hazard rates}
\author{{\bf M. Doostparast}${}^{1,}$\footnote{Corresponding author.\newline
\text{E-mail address:} {\it doustparast@um.ac.ir} (M.
Doostparast). }
\ \ {\bf M. Hashempour}${}^{2}$, and {\bf E. Velayati Moghaddam}${}^{1}$
   \\
{\small {\it $~^{1}$ Department of Statistics, School of
Mathematical Sciences,}}\vspace{-0.2cm}\\ {\small {\it Ferdowsi University of Mashhad, P. O. Box 91775-1159,  Mashhad,\ Iran }}
\vspace{-0.2mm}\\
{\small{\it $~^{2}$ Department of Statistics, School of Science, University of Hormozgan,}}
\vspace{-0.2cm}\\ 
{\small{ \it  Bandar Abbas,\ Iran }}
}

%\author{Majid Hashempour  \footnote{ \newline \noindent Email address: 
%ma.hashempour@hormozgan.ac.ir (M. Hashempour).}
%\vspace{2mm}\\
%{\it\small{ Department of Statistics, School of Science,}}
%\vspace{-0.1cm}\\
%{\it\small{ University of Hormozgan , Bandar Abbas,\ Iran }}}
%

\date{}
\maketitle
\begin{abstract}
In engineering systems, it is usually assumed that lifetimes
of components are independent and identically distributed
(iid). But, the failure of a component results in a higher load on the remaining
components and hence causes the distribution of the surviving components change. For modeling this kind of systems,
the theory of sequential order statistics (SOS) can be used.
Assuming Weibull distribution for lifetimes of components
and {\em conditionally proportional hazard rates} model as a
special case of the SOS theory, the maximum likelihood estimates
of the unknown parameters are obtained in different cases. A new
model, denoted by PTCPHM, as a generalization of the iid case is
proposed, and then statistical inferential methods including point and interval
estimation as well as hypothesis tests under PTCPHM are then developed.
Finally, a real data on failure times of aircraft components, due
to {Mann and  Fertig (1973)}, is analyzed to illustrate the model and inferential methods developed here.
\end{abstract}

\section{Introduction}
A system consisting of $n$ components is said to be a $r$-out-of-$n$
$F$-system if it fails when at least $r$ failures occur (Smith,
2002 and Billinton and Allan, 1992). For $r=1$ and $r=n$, it reduces to series and parallel systems, respectively. Let
$X_1,\cdots,X_n$ be the lifetimes of the components in the system. Then, the
lifetime of a $r$-out-of-$n$ $F$-system is $X_{(r)}$, the
$r$-th order statistic among $X_1,\cdots,X_n$. Thus, order
statistics play an important role in the analysis of these
systems.\par

In the literature, it is usually assumed that the random
variables $X_1,\cdots,X_n$ are independent and identically
distributed (iid). For more details, see Arnold {\em et al.}
(2008). This assumption is violated in many practical
engineering systems since as the components fail sequentially, the
stress on the remaining components would be increased (see,
Cramer and Kamps, 1996). Balakrishnan {\em et al.} (2008) gave
the following example: $\grave{~}\grave{~}$..., the failure of a high-voltage
transmission line will increase the load put on the remaining
high-voltage transmission lines, thus violating the iid
assumption." A method for modeling these systems is through the theory
of sequential order statistics (SOS) (Kamps, 1995). Under this
model, the distribution of the remaining components are changed
when some components fail. {Hashempour and Doostparast (2016) considered Bayesian inference on multiply sequential order statistics from heterogeneous exponential populations with GLR test for homogeneity;} see also Schenk {\em et al}. (2011) considered Bayes estimation and prediction on the basis of multiply Type-II censored data arising from one- and two-parameter exponential distributions;  Shafay et al. (2012) for some additional results in this regard. Since the Weibull distribution is more flexible than the exponential distribution for modeling  failure times as it possesses both
increasing failure rate (IFR) and decreasing failure rate (DFR)
properties, we consider here the problem of estimating the
parameters of the two-parameter Weibull model based on a Type-II
censored sample of SOS.  

The rest of this paper is organized as follows. In Section
2, the likelihood function (LF) associated with SOS coming from a
general class of distribution functions is presented. A new power
trend model, which includes the iid case as a special case, is proposed
in Subsection \ref{label:Description:PTCPHM}. The two-parameter
Weibull distribution is considered in Section 3 in more
detail. More specifically, point estimates as well as
approximate confidence intervals are obtained for the
parameters of the Weibull distribution based on SOS. Assuming Weibull
distribution for the lifetimes of the components, the problem of
hypothesis testing for the new model is discussed in Subsection 3.3.
Also, a test of exponentially for the random variables
$X_1,\cdots,X_n$ against the Weibull model is considered in
Subsection 3.4.
In Section 4, a real data on failure times of aircraft components, due
to Smith (2002), is analyzed to illustrate the model and inferential methods developed here. Finally, some concluding remarks are made in Section 5.

\section{Description of the model and the likelihood function}

In this section, we first describe the conditionally proportional hazard rate model in the setup of SOS. Then, we present the likelihood function for a Type-II censored sample consisting of the first $r$ SOS in the general case, and then its explicit form for the special case of a general exponential family of distributions. Finally, we focus on a special case of the conditionally proportional hazard rate model for SOS called the power trend conditionally proportional hazard model and present the corresponding likelihood function, which is what is used in the subsequent sections to develop inferential methods for the model parameters.

\subsection{The model}

Suppose an engineering system has $n$ components with random lifetimes $X_1,\cdots,X_n$. We assume, to begin with, that these components function independently and have identical lifetimes with cumulative distribution function (cdf) $F_1(x)$, probability density function (pdf) $f_1(x)$, and hazard rate function (hf) $h_1(x)=f_1(x)/(1-F_1(x))$. Next, we assume that, following the first failure, all surviving units face an increased load and continue to function independently but identically distributed with cdf, pdf and hf changed to $F_2(x)$, $f_2(x)$ and $h_2(x)$, respectively, and so on. In general, immediately following the $j$-th failure, the remaining $n-j$ surviving units face an increased load and continue to function independently but distributed with cdf, pdf and hf changed to $F_{j+1}(x)$, $f_{j+1}(x)$ and $h_{j+1}(x)$, respectively, for $j=2,3,\cdots$.\\

In practice, it is quite reasonable to assume that (Balakrishnan {\em et al}., 2008)
\begin{equation}\label{trend:hazard:rates:generally}
h_1(t)<h_2(t)<\cdots<h_n(t).
\end{equation}
One way of considering the model in \eqref{trend:hazard:rates:generally} is through the proportional hazard rate model. Specifically, we take
\begin{equation}\label{h1}
h_j(t)=\alpha_j h_0(t),\ \ \ \mbox{for}\ j=1,2,\cdots,n,
\end{equation}
where $h_0(t)$ is a baseline hazard rate function and  $\alpha_1,\cdots,\alpha_n$ are positive constants. In fact, under this setting, we assume that $F_j(t)=1-[1-F_0(t)]^{\alpha_j}$ for $j=1,\cdots,n$ where $F_0(t)$ is the cdf of the baseline distribution. Thus, under this model, it assumed that when the $j$-th failure occurs, the failure rates of the remaining components in the system are changed from $\alpha_j h_0(t)$ to $\alpha_{j+1} h_0(t)$. In the literature, this model has been termed as {\em conditionally proportional hazard rate model}.

\begin{remark} If we choose the constants $\alpha_j$ in \eqref{h1} such that  $\alpha_1<\alpha_2<\cdots<\alpha_n$, then the {\em restricted sequential order statistics} are obtained.
For more details, see Balakrishnan {\em et al.} (2008) and
the references contained therein.
\end{remark}

\subsection{The likelihood function}
Suppose the first $r$ SOS, denoted by
$\textbf{x}_{*}=(x_1,\cdots,x_r)$, are observed from the model in  \eqref{h1} with pdf and cdf $f(x)$ and $F(x)$, respectively.
Then, the joint pdf of $\textbf{x}_{*}$ is (Cramer and Kamps, 1996)
\begin{eqnarray}
f(x_1,\cdots,x_r)&=&
\frac{n!}{(n-r)!}\left(\prod_{j=1}^{r}\alpha_j\right)\left(\prod_{j=1}^{r-1}[1-F(x_{j})]^{m_j}f(x_{j})\right)\nonumber\\
&&\times\left[1-F(x_{r})\right]^{\alpha_r(n-r+1)-1}f(x_r) ,\label{joint:pdf:sos}
\end{eqnarray}
where ${\bm \alpha}=(\alpha_1,\cdots,\alpha_r)$ and
$m_{j}=(n-j+1)\alpha_{j}-(n-{j})\alpha_{j+1}-1$ for
$j=1,\cdots,r-1$. Let $\mathcal{C}$ be the class of all absolutely continuous distribution functions $F(x;\theta)$ of the form
\begin{eqnarray}
F(x;\theta)&=&1-\exp\left\{-{k_\theta(x)}\right\},~~~~~~~x>0
,\label{cdf:general:class}
\end{eqnarray}
 and hence with pdf $f(x;\theta)=k'_\theta(x)\exp\left\{-{k_\theta(x)}\right\}$, where
$k_\theta(x)$ is the cumulative hazard function, increasing in $x$,
$k'_\theta(x)=\frac{d}{d\theta} k_\theta(x)$, and $\theta$ is a vector of parameters. Substituting \eqref{cdf:general:class} into \eqref{joint:pdf:sos}, the likelihood function (LF) of $\textbf{x}_{*}$ becomes
\begin{eqnarray}
L(\theta;\textbf{x}_{*})
&=&\frac{n!}{(n-r)!}\left(\prod_{j=1}^{r}\alpha_{j}\right)\left(\prod_{j=1}^{r}k'_\theta(x_j)\right)\nonumber\\
&&\times\exp\left[-\left(\sum_{j=1}^{r-1}(m_j+1)k_\theta(x_j)+k_\theta(x_r)\alpha_r(n-r+1)\right)\right].\label{like:general:class}
\end{eqnarray}
The class of distributions $\mathcal{C}$ in \eqref{cdf:general:class} includes
many different lifetime distributions such as the exponential, Weibull and Pareto
models (see AL-Hussaini, 1999). In what follows, we consider
the Weibull distribution under this framework as the working model and develop different inferential methods for the model parameters.  First, we propose a new
alternative model here called the {\it power trend
conditionally proportional hazard model}, denoted by PTCPHM.
\subsection{Description of PTCPHM}\label{label:Description:PTCPHM}
A special case of the conditionally proportional hazard rate
model for SOS in \eqref{h1} is when $\alpha_j=a^j$ for $j=1,\cdots,n$
and $a>0$. We refer to this model as the {\it power trend conditionally
proportional hazard model}, since in this case
\begin{equation}\label{PTCPHM:equation:hazard:rate}
h_{j}(t)=a h_{j-1}(t)=\cdots=a^j h_0(t),\ \ \ \forall\ t>0.
\end{equation}
%Figure \ref{fig:model} shows a prespective of the suggested model.
%\begin{figure}
%\centering
%\includegraphics[angle=0,scale=1]{f4.bmp}
%\includegraphics[width=4in,height=2.5in]{f4.bmp}
%\caption{A prespective for PTCPHM} \label{fig:model}
%\end{figure}
Then, the
LF in \eqref{like:general:class} reduces to
\begin{eqnarray}
L(\theta;\textbf{x}_{*})&=&\frac{n!}{(n-r)!}\left(a^{\frac{r(r+1)}{2}}\right)\left(\prod_{j=1}^{r}k'_\theta(x_j)\right)\nonumber\\
&&\times
\exp\left[-\left(\sum_{j=1}^{r-1}((n-j+1)a^{j}-(n-j)a^{j+1})k_\theta(x_j)+k_\theta(x_r)a^{r}(n-r+1)\right)\right].\nonumber\\
\label{like:general:class:trend}
\end{eqnarray}
For $a>1$, PTCPHM reduces to the restricted SOS, while for
$a=1$, it corresponds to the usual order statistics based on a random sample. For this reason, we will treat the problem of testing the hypothesis $H:
a=1$ against the alternative $K:a>1$ later in Subsection 3.3.

\section{Weibull analysis}
The Weibull distribution has been used extensively in life-testing and reliability studies as it is a flexible lifetime model that includes increasing hazard rate and decreasing hazard rate in addition to the constant hazard rate corresponding to the exponential distribution; see, for example, Johnson {\em et al.} (1994). From the cdf of the Weibull distribution, it is evident that the Weibull distribution is a
member of the class $\mathcal{C}$ in \eqref{cdf:general:class}, with
\begin{equation}\label{k:weibull}
k_\theta(x)=\left(\frac{x}{\sigma}\right)^\beta,
\end{equation}
where $\theta=(\beta,\sigma)$, $\beta$ is the shape parameter and $\sigma$ is the scale parameter. We shall denote this Weibull distribution by $W(\beta,\sigma)$ from now on.
In what follows, we develop inferential methods for the
unknown parameters of the Weibull model based on SOS.
\subsection{Point estimation}
Statistical inference on the basis of SOS has been discussed extensively in the literature; for example, see
Cramer and Kamps (1996, 1998, 2001, 2003), Balakrishnan {\em et al.} (2008), Bedbur (2010), Schenk {\em et al}. (2011), and Shafay {\em et al.} (2012).

Schenk {\em et al}. (2011) considered Bayes
estimation and prediction on the basis of multiply Type-II
censored data arising from one- and two-parameter exponential
distributions; see also Shafay {\em et al.} (2012). Since the Weibull distribution is more flexible than
the exponential distribution for modeling failure times as mentioned above, we discuss here the estimation of parameters of the two-parameter Weibull model based on a Type-II
censored sample of SOS.

From \eqref{like:general:class} and \eqref{k:weibull}, the LF
associated with $\textbf{x}_{*}$ simplifies to
\begin{eqnarray}
L({\bm \alpha},\beta,\sigma;\textbf{x}_{*})&=&\left(\frac{n!}{(n-r)!}\right)\left(\prod_{j=1}^{r}\alpha_j\right)\left(\frac{\beta^r}{\sigma^{r\beta}}\right)\left[\eta(\textbf{x}_{*})\right]^{\beta-1}\nonumber\\
&&\times\exp\left[-\frac{1}{\sigma^\beta}\left(\sum_{j=1}^{r-1}(m_j+1)x_j^\beta+\alpha_r(n-r+1)x_r^\beta\right)\right],\label{pdf:wiebull}
\end{eqnarray}
and the log-likelihood function (LLF) becomes
\begin{eqnarray}
l({\bm \alpha},\beta,\sigma;\textbf{x}_{*})&=&\log \left(\frac{n!}{(n-r)!}\right)+\sum_{j=1}^{r}\log\alpha_{j}+r\log\beta-r\beta\log\sigma+(\beta-1)\log
\eta(\textbf{x}_{*})\nonumber\\
&&-\frac{1}{\sigma^\beta}\left[\sum_{j=1}^{r-1}(m_j+1)x_j^\beta+\alpha_r(n-r+1)x_r^\beta\right],\label{log:like:weibull}
\end{eqnarray}
where $\displaystyle\eta(\textbf{x}_{*}):=\prod_{j=1}^{r}x_j$. We shall now develop the maximum likelihood (ML)
 estimates of the unknown parameters
separately in two cases depending on whether $\alpha$ is known or unknown.
\subsection*{ Case (i): ${\bm \alpha}$ known and $(\beta,\sigma)$ unknown}
Suppose the constants $\alpha_{j}~(1\leq {j}\leq {r})$ in
\eqref{pdf:wiebull} are all known. From \eqref{log:like:weibull}, the
likelihood equations for $\beta$ and $\sigma$ are derived in this case as follows:
\begin{equation}
\left\{
\begin{array}{l}
\frac{r}{\beta}-r\log\sigma+\log
\eta(\textbf{x}_{*})+\frac{\log\sigma}{\sigma^\beta}\left[\sum_{j=1}^{r-1}(m_j+1)x_j^\beta+\alpha_r(n-j+1)x_r^\beta\right]\\
-\frac{1}{\sigma^\beta}\left[\sum_{j=1}^{r-1}(m_j+1)(\log{x_j})(x_j^\beta)+\alpha_r(n-r+1)(\log{x_r})x_r^\beta\right]=0,\\ \\
-\frac{r\beta}{\sigma}+\frac{\beta}{\sigma^{\beta+1}}\left[(\sum_{j=1}^{r-1}(m_j+1)x_j^\beta+\alpha_r(n-r+1)x_r^\beta)\right]=0.
\end{array}
\right.
\end{equation}
After some algebraic calculations not reported here for the sake of brevity, the ML estimate of the scale parameter
$\sigma$ is obtained as
\begin{eqnarray}
\widehat{\sigma}_{1}&=&\left(\frac{(\sum_{j=1}^{r-1}(m_j+1)x_j^{\widehat{\beta}_{1}}+\alpha_r(n-r+1)x_r^{\widehat{\beta}_{1}})}{r}\right)^{\frac{1}{\widehat{\beta}_{1}}},\nonumber
\end{eqnarray}
where $\widehat{\beta}_{1}$ is the ML estimate of the shape parameter $\beta$
obtained numerically by solving the following equation:
\begin{eqnarray}\label{f:h1}
\frac{1}{\beta}=\frac{\sum_{j=1}^{r-1}(m_j+1)(\log
{x_j})(x_j^\beta)+\alpha_r(n-r+1)(\log
{x_r})x_r^\beta}{\sum_{j=1}^{r-1}(m_j+1)x_j^\beta+\alpha_r(n-r+1)x_r^\beta}-\frac{1}{r}\sum_{j=1}^{r}\log
{x_j}.
\end{eqnarray}
For the special case when $r=n$ and $\alpha_{1}=\cdots=\alpha_{n}=1$,
SOS reduces to the usual complete sample. In this case, the ML
estimates of the parameters $\sigma$ and $\beta$ simplify as follows:
\begin{eqnarray*}
\widehat{\sigma}=\sqrt[\widehat{\beta}]{\sum_{j=1}^{n}\frac{x^{\widehat{\beta}}_j}{n}}
\end{eqnarray*}
and $\widehat{\beta}$ as the solution of the equation
\begin{eqnarray}\label{f:h2}
\frac{1}{\widehat{\beta}}+\frac{1}{n}\sum_{j=1}^{n}\log
x_j=\frac{\sum_{j=1}^{n}(\log x_j)x^{\widehat{\beta}}_j}{\sum_{j=1}^{n}x^{\widehat{\beta}}_j},
\end{eqnarray}
respectively, as given in Lehmann and Cassella
(1998, p.468). They proved that Eq. \eqref{f:h2} has a
unique solution.
It should be mentioned here that Balakrishnan and Kateri (2008) extended this result to different forms of censored Weibull data. In an analogous manner, we shall show here that Eq. \eqref{f:h1} has a
unique solution. To end this, let
{\begin{eqnarray}\label{f:h3}
h(\beta)&=&\frac{\sum_{j=1}^{r-1}(m_j+1)(\log
{x_j})(x_j^\beta)+\alpha_r(n-r+1)(\log
{x_r})x_r^\beta}{\sum_{j=1}^{r-1}(m_j+1)x_j^\beta+\alpha_r(n-r+1){x_r}^\beta}-\frac{1}{r}\sum_{j=1}^{r}\log
{x_j}. 
\end{eqnarray}
Then, Eq. \eqref{f:h1} is equivalent to the equation $h(\beta)=1/\beta.$ Notice that \\
\begin{eqnarray}\label{partial}
\frac{\partial
h(\beta)}{\partial{\beta}}
%&=&\left[\sum_{j=1}^{r-1}(m_j+1)(\log
%{x_j})^2 x_j^\beta+\alpha_r(n-r+1)(\log
%{x_r})^2x_r^\beta\right]\nonumber\\
%&&\times\left[\sum_{j=1}^{r-1}(m_j+1)x_j^\beta+\alpha_r(n-r+1){x_r}^\beta\right]\nonumber\\
%&&\times\left[\sum_{j=1}^{r-1}(m_j+1)x_j^\beta+\alpha_r(n-r+1)x_r^\beta\right]^{-2}\nonumber\\
%&&-\left(\frac{\sum_{j=1}^{r-1}(m_j+1)(\log
%{x_j})x_j^\beta+\alpha_r(n-r+1)(\log
%{x_r})x_r^\beta}{\sum_{j=1}^{r-1}(m_j+1)x_j^\beta+\alpha_r(n-r+1){x_r}^\beta}\right)^2+\frac{1}{\beta^2}\nonumber\\
&=&\frac{\sum_{j=1}^{r-1}(m_j+1)(\log
{x_j})^2(x_j^\beta)+\alpha_r(n-r+1)(\log
{x_r})^2x_r^\beta}{\sum_{j=1}^{r-1}(m_j+1)x_j^\beta+\alpha_r(n-r+1){x_r}^\beta}\nonumber\\
&&-\left(\frac{\sum_{j=1}^{r-1}(m_j+1)(\log
{x_j})(x_j^\beta)+\alpha_r(n-r+1)(\log
{x_r})x_r^\beta}{\sum_{j=1}^{r-1}(m_j+1)x_j^\beta+\alpha_r(n-r+1){x_r}^\beta}\right)^2.
\label{derivative:h}  
%&=&Var(\log {Y})+\frac{1}{\beta^2}>0,~~~~\forall\beta>0,\nonumber
\end{eqnarray}  }
{For $a<n/(n-1)$, we prove that Equation \eqref{partial} is positive. To do this, for  $j=1,\cdots,r-1$, let } 
\begin{eqnarray*}
a_j&=&\frac{\sqrt{(m_j+1)}(\log
{x_j})(x_j^{\beta/2})}{\sqrt{\sum_{j=1}^{r-1}(m_j+1)x_j^\beta+\alpha_r(n-r+1){x_r}^\beta}}\\
b_j&=&a_j\\
a_r&=&\frac{\sqrt{\alpha_r(n-r+1)}(\log
{x_r})(x_r^{\beta/2})}{\sqrt{\sum_{j=1}^{r-1}(m_j+1)x_j^\beta+\alpha_r(n-r+1){x_r}^\beta}}\\
b_r&=&a_r.
\end{eqnarray*}
{Therefore, $m_{j}+1=(n-j+1)\alpha_{j}-(n-{j})\alpha_{j+1}>0.$}
Applying the well known Cauchy-Schwartz inequality in \eqref{derivative:h},
%
%where $Y$ is a discrete random variable with probability mass
%function
%\begin{equation}
%P(Y=x_j)=\left\{
%\begin{array}{lcl}
%\frac{(m_j+1) x_{j}^\beta}{\sum_{j=1}^{r-1}(m_j+1)x_{j}^\beta+\alpha_r(n-r+1){x_r}^\beta
%},&~~~&j=1,\cdots,r-1,\\ ~~\\
%\frac{\alpha_r(n-r+1) x_r^\beta}{\sum_{j=1}^{r-1}(m_j+1)x_j^\beta+\alpha_r(n-r+1){x_r}^\beta},&&j=r.
%\end{array}
%\right.
%\end{equation}
we conclude that the function $h(\beta)$ is increasing in $\beta$. Also, it is easy to verify that
\begin{eqnarray}
-\infty&=&\lim_{\beta\to{0}}{h(\beta)}<\frac{1}{r}\sum_{j=1}^{r}\log
{x_{j}}<\log {x_{(r)}}=\lim_{\beta\to{\infty}}{h(\beta)}.\nonumber
\end{eqnarray}
{Hence, the equation $h(\beta)=1/\beta$ has a
unique solution, as required.
 For $a\geq n/(n-1)$, the problem of uniqueness MLE remains open.}

\subsection*{Case (ii): ${\bm \alpha}$ and $(\beta,\sigma)$ are unknown}
In this case, we assume that the vector ${\bm
\alpha}=(\alpha_{1},\cdots,\alpha_{r})$ and the parameters
$\beta$ and $\sigma$ in \eqref{pdf:wiebull} are all unknown. Since
there would be $r+2$ unknown parameters in this case (Cramer and Kamps, 1996), we restrict our attention to the subclass with power trend in
proportionality for the hazard rate function, defined by
$\alpha_{j}=a^{j},~~1\leq{j}\leq{r}$, where $a>1$ is an unknown
parameter. Earlier in Subsection \ref{label:Description:PTCPHM}, we
referred to this subclass as PTCPHM. This assumption implies that
$\alpha_1<\alpha_2<\cdots<\alpha_r$, as supposed by Balakrishnan
{\em et al.} (2008). Thus, the hazard function of lifetimes of
surviving components will increase. This case has been
considered in the literature and is known as {\em order
restricted sequential order statistics}; see, for example,
Balakrishnan {\em et al.} (2008). Under PTCPHM, the LF in
\eqref{pdf:wiebull} becomes
\begin{eqnarray}\label{pdf:wiebull:case2}
L(\beta,\sigma,a;\textbf{x}_{*})&=&\frac{n!}{(n-r)!}\left(a^{\frac{r(r+1)}{2}}\right)\frac{\beta^{r}}{\sigma^{r\beta}}\left[\eta(\textbf{x}_{*})\right]^{\beta-1}\nonumber\\
&&\times\exp\left[-\frac{1}{\sigma^\beta}\left(\sum_{j=1}^{r-1}\bigg[(n-j+1)a^j-(n-j)a^{j+1}\bigg]x_j^\beta+a^{r}(n-r+1)x_r^\beta\right)\right],\nonumber\\
\end{eqnarray}
with the corresponding LLF as
\begin{eqnarray}
l(\beta,\sigma,a;\textbf{x}_{*})&=&\log\left(\frac{n!}{(n-r)!}\right)+\frac{r(r+1)}{2}\log
a+r\log\beta-r\beta
\log\sigma+(\beta-1)\log \eta(\textbf{x}_{*})\nonumber\\
&&-\frac{1}{\sigma^\beta}\left[\sum_{j=1}^{r-1}\left[(n-j+1)a^{j}-(n-j)a^{j+1}\right]x_j^\beta
+a^{r}(n-r+1)x_r^\beta\right].\label{weibull:case2:log:like}
\end{eqnarray}
Thus, the ML estimates of the unknown parameters $\beta$, $\sigma$ and $a$ need to be obtained by solving the likelihood equations
\begin{equation}\label{like:equation:general}
\partial l/\partial\sigma=\partial l/\partial\beta=\partial l/\partial a=0.
\end{equation}
Explicit expressions for the partial derivatives in Eq. \eqref{like:equation:general} are presented in the Appendix.
From Eq. \eqref{like:equation:general} and after some algebraic manipulations, the ML estimate of the
parameter $\sigma$ is
\begin{eqnarray}
\widehat{\sigma}_{2}&=&\left(\frac{
\sum_{j=1}^{r-1}(m_j+1)x_j^{\widehat{\beta}_{2}}+{\widehat{a}_{2}^{r}}(n-r+1)x_r^{\widehat{\beta}_{2}}}{r}\right)^{\frac{1}{\widehat{\beta}_{2}}},\label{wei:case2:mle:sigma}
\end{eqnarray}
where $\widehat{a}_{2}$ and $\widehat{\beta}_{2}$ are the ML estimates of the parameters $a$ and $\beta$, obtained by solving the following equations:
{
\begin{eqnarray}
a=\frac{(r+1)\left(\sum_{j=1}^{r-1}(m_j+1)x_j^{\beta_{2}}+{{a}_{2}^{r}}(n-r+1)x_r^{\beta_{2}}
\right)}{2\left[\sum_{j=1}^{r}(n-j+1)~ja^{j-1}
x_j^{\beta_{2}}-\sum_{j=1}^{r-1}(n-j)(j+1)a^jx_j^{\beta_{2}}\right]}\label{wei:case2:mle:a}
\end{eqnarray}
and
\begin{eqnarray}\label{f:h1:1}
\frac{1}{\beta}=\frac{\sum_{j=1}^{r-1}(m_j+1)(\log
{x_j})(x_j^\beta)+a^r(n-r+1)(\log
{x_r})x_r^\beta}{\sum_{j=1}^{r-1}(m_j+1)x_j^\beta+a^r(n-r+1)x_r^\beta}-\frac{1}{r}\sum_{j=1}^{r}\log
{x_j},
\end{eqnarray}
where $m_{j}=(n-j+1)a^{j}-(n-{j})a^{j+1}-1$ for
$j=1,\cdots,r-1$.}

\subsection{Approximate interval estimation}
Since Eq. \eqref{f:h1} has a unique solution, the ML
estimates are consistent, asymptotically normal, and efficient.
Therefore, from Lehmann and Casella (1998, p. 463, Theorem 5.1),
the random vector
$\sqrt{n}(\hat{\beta}_2-\beta,\hat{\sigma}_2-\sigma,\hat{a}_2-a)^T$
converges to the multivariate normal
$\textbf{N}_3(\textbf{0}^T,\left[\textbf{I}(\beta,\sigma,a)\right]^{-1})$
as $n$ goes to infinity, where $\textbf{0}^T=(0,0,0)$ and
$\textbf{I}(\beta,\sigma,a)$ is the Fisher information matrix given by
\begin{equation}\label{fisher:matrix}
\textbf{I}(\beta,\sigma,a)=\left[
\begin{array}{ccc}
w_{11}&w_{12}&w_{13}\\
w_{21}&w_{22}&w_{23}\\
w_{31}&w_{32}&w_{33}\\
\end{array}
\right],
\end{equation}
where $w_{11}=-E\left[\partial^2\log L/\partial\beta^2\right]$,
$w_{12}=w_{21}=-E\left[\partial^2\log
L/\partial\beta\partial\sigma\right]$,
$w_{13}=w_{31}=-E\left[\partial^2\log L/\partial\beta\partial
a\right]$, $w_{22}=-E\left[\partial^2\log
L/\partial\sigma^2\right]$, $w_{23}=w_{32}=-E\left[\partial^2\log
L/\partial\sigma\partial a\right]$, and
$w_{33}=-E\left[\partial^2\log L/\partial a^2\right]$. Explicit
expressions for $w_{ij}$, $i,j=1,2,3$, are presented in the Appendix.
%One can easily verify that
%\begin{equation}\label{inverse:fisher:matrix}
%\left[\widehat\textbf{I}(\beta,\sigma,a)\right]^{-1}=\frac{1}{|\widehat{\textbf{I}}(\beta,\sigma,a)|}\left[
%\begin{array}{ccc}
%b_{11}&b_{12}&b_{13}\\
%b_{21}&b_{22}&b_{23}\\
%b_{31}&b_{32}&b_{33}
%\end{array}
%\right],
%\end{equation}
%where
%\begin{equation}\label{wei:case2:elements:inverse:matrix}
%\left\{
%\begin{array}{l}
%|{\textbf{I}}(\beta,\sigma,a)|=w_{11}(w_{22}w_{33}-w_{23}w_{32})-w_{12}(w_{21}w_{33}-w_{23}
%w_{31})+w_{13}(w_{21}w_{32}-w_{22}w_{31}),\\
%b_{11}=w_{22}w_{33}-w_{23}w_{32},\\
%b_{12}=b_{21}=-(w_{12}w_{33}-w_{32}w_{13}),\\
%b_{13}=b_{31}=w_{12}w_{23}-w_{22}w_{13},\\
%b_{22}=w_{11}w_{33}-w_{13}w_{31},\\
%b_{23}=b_{32}=-(w_{11}w_{23}-w_{13}w_{21}),\\
%b_{33}=w_{11}w_{22}-w_{12}w_{21}.
%\end{array}
%\right.
%\end{equation}
Since it is not possible to obtain the expectations involved in $w_{ij}$, we will use the observed Fisher information,  denoted
by $\hat{I}(\beta,\sigma,a)$,  obtained
by replacing $\beta$, $\sigma$ and
$a$ by the corresponding ML estimates based on SOS into
\eqref{fisher:matrix}.
%More specifically, we
%have {$\hat{I}(\beta,\sigma,a)=[[\hat{w}_{ij}]]$} where
%\[\hat{w}_{11}=-{\partial^2\log L}/{\partial\beta^2}
%|_{\beta=\hat{\beta}_{2},\sigma=\hat{\sigma}_{2},a=\hat{a}_{2}},\]
%\[\hat{w}_{12}=\hat{w}_{21}=-{\partial^2\log L
%}/{\partial\beta\partial \sigma}
%|_{\beta=\hat{\beta}_{2},\sigma=\hat{\sigma}_{2},a=\hat{a}_{2}},\]
%\[\hat{w}_{13}=\hat{w}_{31}=-{\partial^2\log
%L}/{\partial\beta\partial
%a}|_{\beta=\hat{\beta}_{2},\sigma=\hat{\sigma}_{2},a=\hat{a}_{2}},\]
%\[\hat{w}_{22}=-{\partial^2\log L}/{\partial\sigma^2}
%|_{\beta=\hat{\beta}_{2},\sigma=\hat{\sigma}_{2},a=\hat{a}_{2}},\]
%\[\hat{w}_{23}=\hat{w}_{32}=-{\partial^2\log
%L}/{\partial\sigma\partial
%a}|_{\beta=\hat{\beta}_{2},\sigma=\hat{\sigma}_{2},a=\hat{a}_{2}}\]
%\[\hat{w}_{33}=-{\partial^2\log L}/{\partial a^2}
%|_{\beta=\hat{\beta}_{2},\sigma=\hat{\sigma}_{2},a=\hat{a}_{2}}.\]
Hence, the approximate $100(1-\gamma)$\% equi-tailed
confidence intervals for $\beta$, $\sigma$ and $a$ are,
respectively, given by
\begin{equation}\label{approx:inter:beta}
\left(\hat{\beta}-z_{1-\gamma/2}\sqrt{\frac{1}{|\widehat{\textbf{I}}(\beta,\sigma,a)|}\hat{b}_{11}},
\hat{\beta}+z_{1-\gamma/2}\sqrt{\frac{1}{|\widehat{\textbf{I}}(\beta,\sigma,a)|}\hat{b}_{11}}\right),
\end{equation}
\begin{equation}\label{approx:inter:sigma}
\left(\hat{\sigma}-z_{1-\gamma/2}\sqrt{\frac{1}{\left|\widehat{\textbf{I}}(\beta,\sigma,a)\right|}\hat{b}_{22}},
\hat{\sigma}+z_{1-\gamma/2}\sqrt{\frac{1}{|\widehat{\textbf{I}}(\beta,\sigma,a)|}\hat{b}_{22}}\right),
\end{equation}
and
\begin{equation}\label{approx:inter:a}
\left(\hat{a}-z_{1-\gamma/2}\sqrt{\frac{1}{|\widehat{\textbf{I}}(\beta,\sigma,a)|}\hat{b}_{33}},
\hat{a}+z_{1-\gamma/2}\sqrt{\frac{1}{|\widehat{\textbf{I}}(\beta,\sigma,a)|}\hat{b}_{33}}\right),
\end{equation}
where $z_\gamma$ is the $\gamma$-percentile of the standard
normal distribution,
$b_{11}=w_{22}w_{33}-w_{23}w_{32}$, $b_{12}=b_{21}=-(w_{12}w_{33}-w_{32}w_{13})$, $b_{13}=b_{31}=w_{12}w_{23}-w_{22}w_{13}$,
$b_{22}=w_{11}w_{33}-w_{13}w_{31}$, $b_{23}=b_{32}=-(w_{11}w_{23}-w_{13}w_{21})$, $b_{33}=w_{11}w_{22}-w_{12}w_{21}$, and
 $\hat{b}_{ij}$ is obtained by replacing $w_{ij}$
by the corresponding ML estimate, denoted by $\hat{w}_{ij}$. For simultaneous confidence intervals, there
are various methods. For example, the Bonferroni simultaneous
confidence intervals for the unknown parameters are obtained from
\eqref{approx:inter:beta}, \eqref{approx:inter:sigma} and
\eqref{approx:inter:a} by replacing $z_{1-\gamma/2}$ by
$z_{1-\gamma/6}$. One can also use the recent approach of Casella and
Hwang (2012) to obtain more accurate simultaneous confidence sets.

\begin{remark}
Let $g(\beta,\sigma,a)$ be an arbitrary measurable function of the three
parameters $\beta,\sigma$ and $a$. By the use of multivariate version of
delta method, we obtain
\begin{equation}\label{asy:cdf:function}
\sqrt{n}\left(g(\hat{\beta}_{2},\hat{\sigma}_{2},\hat{a}_{2})-g(\beta,\sigma,a)
\right)\stackrel{D}{\to}\textbf{N}(\textbf{0}^T,\nabla
g(\beta,\sigma,a) \left[\textbf{I}(\beta,\sigma,a)\right]^{-1}\nabla
g(\beta,\sigma,c)^T)
\end{equation}
as $n\to \infty$, where $\nabla g(\beta,\sigma,a)$ denotes the gradient
of the function $g(\beta,\sigma,a)$.
\end{remark}
As an example, suppose we want to estimate the baseline
survival function of the lifetime at a fixed time (say, $t_0$), i.e.,
$S(t_0;\beta,\sigma)=\exp\left\{-({t_0}/{\sigma})^\beta\right\}$.
Then, from \eqref{asy:cdf:function}, we conclude that, as $n\to\infty$,
$\sqrt{n}(S(t_0;\hat{\beta}_{2},\hat{\sigma}_{2})-S(t_0;\beta,\sigma))$
tends to a normal distribution with mean zero and variance
\[\frac{e^{-2({t_0}/{\sigma})^\beta}\left({t_0}/{\sigma}\right)^{2\beta}}{|\widehat{\textbf{I}}(\beta,\sigma,a)|}\left[b_{11}\log^2\left(\frac{t_0}{\sigma}\right)-2b_{21}\left(\frac{\beta}{\sigma}\right)
\log\left(\frac{t_0}{\sigma}\right)+b_{22}\left(\frac{\beta}{\sigma}\right)^2\right],\]
which readily yields an approximate $100(1-\gamma)\%$ equi-tailed
confidence interval for $S(t_0;\beta,\sigma)$ as
\[
S(t_0;\hat{\beta}_{2} ,\hat{\sigma}_{2})\pm z_{1-\gamma/2}\sqrt{\frac{e^{-2({t_0}/{\sigma})^\beta}\left({t_0}/{\sigma}\right)^{2\beta}}{|\widehat{\textbf{I}}(\beta,\sigma,a)|}\left[\hat{b}_{11}\log^2\left(\frac{t_0}{\sigma}\right)-2\hat{b}_{21}\left(\frac{\beta}{\sigma}\right)
\log\left(\frac{t_0}{\sigma}\right)+\hat{b}_{22}\left(\frac{\beta}{\sigma}\right)^2\right]}.\]

\subsection{Hypothesis testing for PTCPHM}
In the preceding section, we assumed a power trend for
$\alpha_{j}$ as $\alpha_{j}=a^{j}$ $(1\leq{j}\leq{r})$. Here, we assume that $a\geq 1.$
Under the null hypothesis $H:a=1$, we assume that the failure of a component
does not effect the distribution of the remaining components,
while under the alternative $K:a>1$, we have an increase in the hazard rates of the
remaining components. Therefore, the problem of testing the null
hypothesis $H:a=1$ against the alternative $K:a>1$ is of natural interest. From Lehmann and Romano (2005, p. 513), a generalized likelihood ratio (GLR) test
has the rejection region as $\left\{\textbf{x}_{*}: \Lambda < k\right\}$,
where
\begin{eqnarray*}
\Lambda&=&\frac{\sup_{\beta>0,\sigma>0} L(\beta,\sigma,1;\textbf{x}_{*})}{\sup _{\beta>0,\sigma>0,a\geq 1}L(\beta,\sigma,a;\textbf{x}_\star)}\\
&=&\left(\frac{\widehat{\beta}_{1 }}{\widehat{\beta}_{2}}\right)^{r}\left(\frac{\sum_{j=1}^{r-1}\big[(n-j)\hat{a}_2^j-(n-j+1)\hat{a}_2^{j+1}\big]x_j^{\widehat{\beta}_{2}}+{\widehat{a}_{2}}^{r}(n-r+1)x_r^{\widehat{\beta}_{2}}}{\sum_{j=1}^{r-1}x_j^{\widehat{\beta}_{1}}+(n-r+1)x_r^{\widehat{\beta}_{1}}}\right)^{r}\left({\widehat{a}_{2}}^{-\left[\frac{r(r+1)}{2}\right]}\right)\nonumber\\
&&\times\left[\eta(\textbf{x}_\star)\right]^{\left(\widehat{\beta}_{1}-\widehat{\beta}_{2}\right)}.
\end{eqnarray*}
Under the null hypothesis $H$ and the usual regularity conditions (see Lehmann and Cassella, 1998), $-2\log \Lambda$ has asymptotically the chi-square distribution with $1$ degree of freedom. Thus, for large $n$, the rejection region of the GLR test of size $\gamma$ is
\begin{eqnarray}\label{rejection:GLR:large:n}
-2\log \Lambda > \chi^{2}_{1,1-\gamma},
\end{eqnarray}
where $\chi^2_{\upsilon,\gamma}$ is the $\gamma$-th precentile of the chi-square distribution with $\upsilon$ degrees of freedom. Also, the actual level of the GLR test may be obtained by means of a Monte Carlo
(MC) simulation study for given $a,\beta$ and $\sigma$.

\subsection{ Exponential baseline distribution}
When $\beta=1$, the $W(\beta,\sigma)$ -distribution reduces to
the exponential distribution, denoted by $Exp(\sigma)$. Now, suppose
 the baseline distribution is $Exp(\sigma)$. Then,
the LF in \eqref{pdf:wiebull} simplifies to
\begin{eqnarray}\label{pdf;exponential}
L(\sigma,\bm \alpha;\textbf{x}_{*})&=&\frac{n!}{(n-r)!}\left(\prod_{j=1}^{r}\alpha_{j}\right)\sigma^{-r}\nonumber\\
&&\times\exp\left[-\frac{1}{\sigma}\left(\sum_{j=1}^{r}(m_j+1)x_j+\alpha_{r}(n-r+1)x_r\right)\right]\label{like:exponential}.
\end{eqnarray}
\noindent For $\bm \alpha$ known, the ML estimate of $\sigma$ is
\begin{eqnarray*}
\widehat{\sigma}_{1,E}&=&\frac{\sum_{j=1}^{r}(m_j+1)x_j+\alpha_r(n-r+1)x_r}{r}.
\end{eqnarray*}
In the case when $\bm \alpha$ is unknown and $\alpha_j=a^{j}$ $(1\leq{j}\leq{r})$, we have
\begin{eqnarray}\label{pdf;exponential;1}
L(\sigma,a;\textbf{x}_{*})&=&\frac{n!}{(n-r)!}a^{\frac{r(r+1)}{2}}\sigma^{-r}\nonumber\\
&&\times\exp\left[-\frac{1}{\sigma}\left(\sum_{j=1}^{r}(m_j+1)x_j+a^{r}(n-r+1)x_r\right)\right].
\end{eqnarray}
{Thus,  the ML estimates of $\sigma$ and $a$ based on SOS are derived by solving the following equations:
\begin{eqnarray*}
\sigma&=&\frac{\sum_{j=1}^{r}(n-j+1){a}^{j}x_j-\sum_{j=1}^{r-1}(n-j){a}^{j+1}x_j}{r}
\end{eqnarray*}
and
\begin{eqnarray*}
a&=&\frac{r(r+1){\sigma}}{2\left[\sum_{j=1}^{r}(n-j+1)(j)a^{j-1}x_j-\sum_{j=1}^{r-1}(n-j)(j+1)a^{j}x_j\right]}\\
&=&\frac{r(r+1)}{2},
\end{eqnarray*} }
respectively. Similarly, a GLR test of size $\gamma$ for the null hypothesis $H: a=1$
against the alternative $K: a>1$ has its critical region as
\begin{eqnarray}\label{GLR:wei:expo:rejection}
-2\log \Lambda > \chi^{2}_{1,\gamma},
\end{eqnarray}
where
\begin{eqnarray*}
\Lambda&=&\left(\frac{\sum_{j=1}^{r-1}(m_j+1)x_j+{\widehat{a}_{2,E}}^{r}(n-r+1)x_r}
{\sum_{j=1}^{r-1}x_j+(n-r+1)x_r}\right)^{r}\left({\widehat{a}_{2,E}}^{-\left[\frac{r(r+1)}{2}\right]}\right).
\end{eqnarray*}
\begin{remark}
Schenk {\it et al}. (2011) and Shafay {\em et al.} (2012) discussed Bayesian estimation and prediction based on sequential order statistics from the exponential distribution.
\end{remark}

\section{Aircraft data set}
To demonstrate the performance of the results obtained in Section
3, we present an illustrative example in this section.\\
Smith (2002, p. 130) gave failure times of aircraft components for a life-test, originally due to Mann and Fertig (1973). In the test, $n=13$ components were placed in a Type-II censored life test in which the failure times of first 10 components to fail were observed (in hours) as
$$0.22, 0.50, 0.88, 1.00, 1.32, 1.33, 1.54, 1.76, 2.50, 3.00.$$
Assuming that the lifetimes of the components are iid with an exponential distribution (i.e., $\beta=1$ and $a=1$  in { Eq. \eqref{pdf:wiebull:case2}),}  the ML estimate of the mean is obtained to be $\hat\sigma=2.305$.

\begin{table}
\caption{Fitted models for failure time of aircraft components}
\begin{tabular}{lccccc}
  \hline
Model  &  & ML estimates &  & Log-likelihood & AIC\\
&$\hat{\beta}$ & $\hat{\sigma}$ & $\hat{a}$\\
&&&&& \\
\hline
$ a=1,  \beta=1$ & - & 2.3050 & - & $-18.3508$&  38.7016\\
$ \beta=1$ & - & 2.9704 & 1.04936 &  $-18.2372$& 40.4743\\
$ a=1$ & 1.41746  & 2.27315 & - & $-17.6335$ & 39.2670\\
$   a\neq 1, \beta>0, \sigma>0$ & 2.02392 & 1.25749 & 0.823473 &  $-16.7801$&39.5602\\    \hline
\end{tabular}
\end{table}

We analyzed these data with different models, and the results of the fitted models are presented in Table 1. Adopting Akaike information criterion (AIC), we conclude that the exponential
model is still supported for the failure times. But, a PTCPHM is also found to be suitable. The inverse of observed Fisher information is
obtained from \eqref{fisher:matrix} as
\begin{equation}\label{example:inverse:fisher:matrix}
I^{-1}=\left[
\begin{array}{ccc}
0.520823& -0.155695& -0.0674688\\
-0.155695& 0.16624&  0.0404666\\
-0.0674688& 0.0404666& 0.0139039
\end{array}
\right].
\end{equation}
Thus, the approximate $95\%$ confidence intervals for the
unknown parameters $\beta$, $\sigma$ and $a$ in
\eqref{log:like:weibull} are obtained from
\eqref{approx:inter:beta}, \eqref{approx:inter:sigma},
\eqref{approx:inter:a} and \eqref{example:inverse:fisher:matrix}
to be $(0.609424,3.43842)$, $(0.458347, 2.05663)$ and
$(0.592359,1.05459)$, respectively. We also obtain an
approximate simultaneous $95\%$ confidence region for the three
parameters from \eqref{approx:inter:beta},
\eqref{approx:inter:sigma} and \eqref{approx:inter:a} by
replacing $z_{1-(0.05/2)}$ by $z_{1-(0.05/6)}=2.39398$ to be
\[(\beta,\sigma,a)\in(0.299101,3.74874)\times (0.283025,2.23196)\times (0.541656,1.10529).\]

\section{Concluding remarks}
In this paper, by considering an engineering system with $n$ components, we have considered the situation when failures of components increase the load on surviving components thus changing their lifetime distribution. For modeling this situation, we have described the conditionally proportional hazard rates model in the framework of sequential order statistics. We have then focused  on a special case of this model called the power trend conditionally proportional hazard model, and developed inferential methods based on a Type-II censored data from this model under the assumption of Weibull distributed lifetimes. Using a well-known data of { Mann and Fertig (1973)}  on failure times of aircraft components, we have illustrated the developed inferential results. Since the framework presented here is applicable for the general exponential family of distributions of the form in \eqref{cdf:general:class}, inferential results analogous to those for the Weibull distribution here can be developed for some other lifetime distributions of interest such as Pareto. It will also be of interest to develop point and interval prediction methods under this general framework. Work on these problems is currently under progress and we hope to report these findings in a future paper.
\section*{Acknowledgements}
The authors are grateful to anonymous referees and the Associate Editor for their useful suggestions and comments on
an earlier version of this paper, which resulted in this improved version of the manuscript.
\section*{References}
\begin{description}
\item AL-Hussaini, E. K. (1999) Predicting observable from a general class of distributions,
{\it Journal of Statistical Planning and Inference}, {\bf79}, 79-91.
\item Arnold, B. C., Balakrishnan, N. and Nagaraja, H. N. (2008) {\it A First Course in Order Statistics}, Clasic Edition, SIAM, Philadelphia.
\item Balakrishnan, N., Beutner, E. and Kamps, U. (2008) Order restricted inference for sequential $k$-out-of $n$ systems, {\it Journal of Multivariate Analysis,} {\bf 99}, 1489-1502.
\item Balakrishnan, N. and Kateri, M. (2008) On the maximum likelihood estimation of parameters of Weibull distribution based on complete and censored data, {\it Statistics \& Probability Letters}, {\bf 78}, 2971-2975.
\item Bedbur, S. (2010) UMPU test based on sequential order statistics,
{\it Journal of Statistical Planning and Inference,} {\bf 140},
2520-2530.
\item Billinton, R. and Allan, R. (1992) {\it Reliability of Engineering Systems: Concepts and Techniques}, Second edition, Springer-Verlag, New York.
\item Casella, C. and Hwang, J. T. G. (2012) Shrinkage confidence procedures,  {\it Statistical Science,} {\bf 27}, 51-60.
\item Cramer, E. and Kamps, U. (1996) Sequential order statistics and $k$-out-of-$n$ systems
with sequentially adjusted failure rates, {\it Annals of the
Institute of Statistical Mathematics} {\bf 48}, 535-549.
\item Cramer, E. and Kamps, U. (1998) Sequential $k$-out-of-$n$ systems with Weibull
components, {\it Economic Quality Control}, {\bf 13}, 227-239.
\item Cramer, E. and Kamps, U. (2001) Sequential $k$-out-of-$n$ systems, In: N. Balakrishnan
and C.R. Rao (Eds.), {\it Handbook of Statistics}, Vol. {\bf 20}, {\it Advances
in Reliability}, pp. 301-372, North-Holland, Amsterdam.
\item Cramer, E. and Kamps, U. (2003) Marginal distributions of sequential and
generalized order statistics, {\it Metrika}, {\bf 58}, 293-310.
\item Johnson, N. L., Kotz, S. and Balakrishnan, N. (1994) {\it Continuous Univariate Distributions}-Vol. 1, Second edition, John Wiley \& Sons, New York.

{ \item Hashempour, M and Doostparast, M. (2016) considered Bayesian inference on multiply sequential order statistics from heterogeneous exponential populations with GLR test for homogeneity, {\it Communications in Statistics-Theory and Methods},  Doi.10.1080/03610926.2016.1175625, 2016.}

\item Kamps, U. (1995) {\it A Concept of Generalized Order Statistics,} Teubner, Stuttgart, Germany.
\item  Lehmann, E. L. and Casella, G. (1998) {\it Theory of Point Estimation}, Second edition, Springer-Verlag, New York.
\item  Lehmann, E. L. and Romano, J. P. (2005) {\it Testing Statistical Hypothesis}, Third edition, Springer-Verlag, New York.
 {\item Mann, N. R. and Fertig, K. W. (1973) Tables for obtaining Weibull confidence bounds and tolerance bounds based on best linear invariant estimates of parameters of the extreme value distribution, {\it Technometrics}, {\bf 15}, 87-101.}
\item Schenk, N. Burkschat, M., Cramer E. and Kamps, U. (2011) Bayesian estimation and prediction with multiply Type-II censored
samples of sequential order statistics from one- and two-parameter
exponential distributions, {\it Journal of Statistical Planning
and Inference,} {\bf 141}, 1575-1587.
\item Shafay, A. R., Balakrishnan, N. and Sultan, K. S. (2012) Two-sample Bayesian prediction for sequential order statistics from exponential distribution based on multiply Type-II censored samples, {\it Journal of Statistical Computation and Simulation} (to appear).
\item Smith, P. J. (2002). {\it Analysis of Failure and Survival Data}, Chapman \& Hall/CRC Press, Boca Raton, Florida.

\end{description}
\section*{Appendix}
From Eq. \eqref{weibull:case2:log:like}, explicit expressions for the partial derivatives in Eq. \eqref{like:equation:general} and for $w_{ij}$, $i,j,=1,2,3$, in \eqref{fisher:matrix} are obtained as follows:
\begin{eqnarray*}
\frac{\partial l}{\partial \sigma}&=&-\frac{r\beta}{\sigma}+\frac{\beta}{\sigma^{\beta+1}}\left[\sum_{j=1}^{r-1}
\bigg\{(n-j+1)a^j-(n-j)a^{j+1}\bigg\}x_j^\beta+a^r(n-r+1)x_r^\beta\right],\label{wei:h1}\\
\frac{\partial l}{\partial \beta}&=&\frac{r}{\beta}-r\log\sigma+\log\eta({\bf x_\star})+\frac{\log\sigma}{\sigma^\beta}\left[\sum_{j=1}^{r-1}
\bigg\{(n-j+1)a^j-(n-j)a^{j+1}\bigg\}x_j^\beta+a^r(n-j+1)x_r^\beta\right]\nonumber\\
&&-\frac{1}{\sigma^\beta}\left[\sum_{j=1}^{r-1}\bigg\{(n-j+1)a^j-(n-j)a^{j+1}\bigg\}
(\log{x_j})x_j^\beta+a^r(n-r+1)(\log{x_r})x_r^\beta\right],\label{wei:h2}\\
\frac{\partial l}{\partial
a}&=&\frac{r(r+1)}{2a}-\frac{1}{\sigma^\beta}\left[\sum_{j=1}^{r-1}\bigg\{(n-j+1)~ja^{j-1}-(n-j)(j+1)a^{j}\bigg\}
x_j^\beta+(n-r+1)~ra^{r-1}x_r^\beta\right];\\
w_{11}&=&-\frac{\partial^{2}l}{\partial{\beta^2}}\\
&=&\frac{r}{\beta^2}+\frac{(\log\sigma)^2}{\sigma^\beta}\left[\sum_{j=1}^{r-1}(m_j+1)x_{j}^\beta+a^r(n-r+1)x_r^\beta\right]\\
&&-\frac{2\log\sigma}{\sigma^\beta}\left[\sum_{j=1}^{r-1}(m_j+1)(\log{x_j})x_{j}^\beta+a^r(n-r+1)(\log{x_r})x_r^\beta\right]\\
&&+\frac{1}{\sigma^\beta}\left[\sum_{j=1}^{r-1}(m_j+1)(\log{x_j})^{2}x_{j}^\beta+a^r(n-r+1)(\log{x_r})^{2}x_r^\beta\right],
\end{eqnarray*}
\begin{eqnarray*}
w_{12}&=&-\frac{\partial^{2}l}{\partial\beta\partial\sigma}\\
&=&\frac{r}{\sigma}+\frac{\beta\log\sigma-1}{\sigma^{\beta+1}}\left[\sum_{j=1}^{r-1}(m_j+1)x_{j}^\beta+a^r(n-r+1)x_r^\beta\right]\\
&&-\frac{\beta}{\sigma^{\beta+1}}\left[\sum_{j=1}^{r-1}(m_j+1)(\log{x_j})x_{j}^\beta+a^r(n-r+1)(\log{x_r})x_r^\beta\right],
\end{eqnarray*}
\begin{eqnarray*}
w_{13}&=&-\frac{\partial^{2}l}{\partial\beta\partial{a}}\\
&=&-\frac{\log\sigma}{\sigma^\beta}\left[\sum_{j=1}^{r-1}\bigg\{(n-j+1)ja^{j-1}-(n-j)(j+1)a^{j}\bigg\}x_j^\beta+(n-r+1)ra^{r-1}x_r^\beta\right]\\
&&+\frac{1}{\sigma^\beta}\bigg[\sum_{j=1}^{r-1}\bigg\{(n-j+1)ja^{j-1}-(n-j)(j+1)a^{j}\bigg\}(\log{x_j})x_j^\beta\\
&&+(n-r+1)ra^{r-1}(\log{x_r})x_r^\beta\bigg],
\end{eqnarray*}
\begin{eqnarray*}
w_{22}&=&-\frac{\partial^{2}l}{\partial{\sigma}^2}\\
&=&\frac{-r\beta}{\sigma^2}+\beta(\beta+1)\sigma^{-(\beta+2)}\left[\sum_{j=1}^{r-1}(m_j+1)x_{j}^\beta+a^r(n-r+1)
x_r^\beta\right],
\end{eqnarray*}
\begin{eqnarray*}
w_{23}&=&-\frac{\partial^{2}l}{\partial\sigma\partial{a}}\\
&=&-\frac{\beta}{\sigma^{\beta+1}}\bigg[\sum_{j=1}^{r-1}\bigg\{(n-j+1)ja^{j-1}-(n-j)(j+1)a^{j}\bigg\}(\log{x_j})x_j^\beta\\
&&+(n-r+1)ra^{r-1}(\log{x_r})x_r^\beta\bigg],
\end{eqnarray*}
\begin{eqnarray*}
w_{33}&=&-\frac{\partial^{2}l}{\partial a ^2}\\
&=&\frac{r(r+1)}{2 a^2}+\sigma^{-\beta}\bigg[\sum_{j=1}^{r-1}\bigg\{(n-j+1)j(j-1) a^{j-2}-(n-j)(j+1)j a^{j-1}\bigg\}x_j^\beta\\
&& +(n-r+1)r(r-1) a^{r-2} x_r^\beta\bigg].
\end{eqnarray*}

\end{document}